\documentclass[a4paper,12pt,leqno]{article}
\usepackage{amsmath,amsfonts}
\usepackage{graphicx}

\newtheorem{theorem}{Theorem}[section]
\newtheorem{proposition}[theorem]{Proposition}
\newtheorem{lemma}[theorem]{Lemma}

\newtheorem{definition}{Definition}[section]

\newcommand{\T}{\mathbb{T}}
\newcommand{\Z}{\mathbb{Z}}
\newcommand{\R}{\mathbb{R}}

\newcommand{\A}{\mathcal{A}}
\newcommand{\B}{\mathcal{B}}

\newcommand{\D}{\mathcal{D}}
\newcommand{\F}{\mathcal{F}}
\newcommand{\G}{\mathcal{G}}

\begin{document}

\begin{center}
{\large\bfseries
Counting Berg partitions}
\end{center}
\begin{center}
{\bf Artur Siemaszko,  Maciej P. Wojtkowski}
\end{center}
\begin{center}
{\it Department of Mathematics and Informatics,
University of Warmia and Mazury,}
\end{center}
\begin{center}
{\it ul. Zolnierska 14,
10-561 Olsztyn, Poland}
\end{center}
\begin{center}
{\it E-mail addresses}:
{\tt artur@uwm.edu.pl, \tt wojtkowski@matman.uwm.edu.pl}
\end{center}

\begin{center}
22 July 2010
\end{center}

%{\it MSC:} 37D;70D;

\begin{abstract}
We call a Markov partition of a two dimensional hyperbolic toral automorphism
a {\it Berg partition} if it contains just two rectangles. We describe all Berg
partitions for a given hyperbolic toral automorphism. In particular there are
exactly $(k+n+l+m)/2$ nonequivalent Berg partitions with the same connectivity
matrix $(k,l,m,n)$.
\end{abstract}

\section{Introduction}
Markov partitions play an important role in the theory
of dynamical systems. They were introduced by Adler and Weiss
in the seminal paper \cite{A-W}, in the context of toral automorphisms.
The notion was further developed by Sinai, \cite{Si1}, \cite{Si2} and Bowen, \cite{B},
becoming a principal scenario for deterministic systems with stochastic
behavior.

In this paper we go back to the original setting of Adler and Weiss,
of a two-dimensional toral automorphism,
and consider Markov partitions with two geometric rectangles.
Such partitions are not generating, but they can be routinely refined to
a generating partition, and there are distinctive advantages to
having partitions of just two elements.

We propose to call such partititions {\it Berg partitions} based on
the following historical comment by Roy Adler
{\it(Kenneth) ``Berg in his Ph.D.
thesis [Be] was the first to discover a Markov partition of a smooth
domain under the action of a smooth invertible map: namely, he constructed
Markov partitions for hyperbolic automorphisms acting on the two
dimensional torus.''
\hfill R. Adler, BAMS, 1998 }, \cite{Ad}.

Our goal is to find and classify all such partitions for a given
toral automorphism $\D$. The Berg partitions differ by shapes and their
placement in the torus. The shapes are shared by all hyperbolic
automorphisms in the centralizer of $\D$ in $GL(2,\mathbb Z)$.
Their number $N$ is related to  the period of the continued fraction
expansion of the slope of the eigenvector of $\D$.
Each Berg partition comes with
the connectivity matrix $C$ with nonnegative entries and $det C = \pm 1$.
There are exactly $N$ different connectivity matrices for a given
automorphism $\D$.

We call two Berg partitions equivalent if there is a
continuous toral map
commuting with $\D$ which takes one into the other.
A Berg partition can be translated into another,
nonequivalent Berg partition. We prove (Theorem 4)
that the number of such translations is equal to one half the sum
of entries of the connectivity matrix.

In the last Section we describe symmetries of Berg partitions.
They are present for reversible toral automorphisms.
The full symmetry group of a toral automorphism was studied by Baake and Roberts \cite{B-R}.

Let us describe the results  that proceeded our work.
Snavely studied in \cite{Sn} the connectivity matrices of Markov
partitions for hyperbolic automorphisms of $\T^2$.
He found that for Berg partitions the connectivity matrices
are conjugated to the dynamics. He also found a way to list all
such matrices and hence to classify the shapes of Berg partitions.
He relied on the result of Adler, \cite{Ad}, that such partitions
are indeed present for any toral automorphism. Manning gave
a powerful generalization of this to $\T^n$, \cite{M}.

After we presented our work at the November 2009 conference
Progress in Dynamics at IHP in Paris
we learned from Pascal Hubert about the work of
Anosov, Klimenko and Kolutsky \cite{A-K-K} which pursues similar goals.
While our paper has intersections with all the previous works,
we give an independent presentation.

\section{Bi-partitions of the torus}\label{bipart}

We reserve the term {\it rectangle} to rectangles $R \subset \mathbb R^2$
\begin{figure}[h]{Fig. 1}
  \centering
    \includegraphics[width=0.7\textwidth]{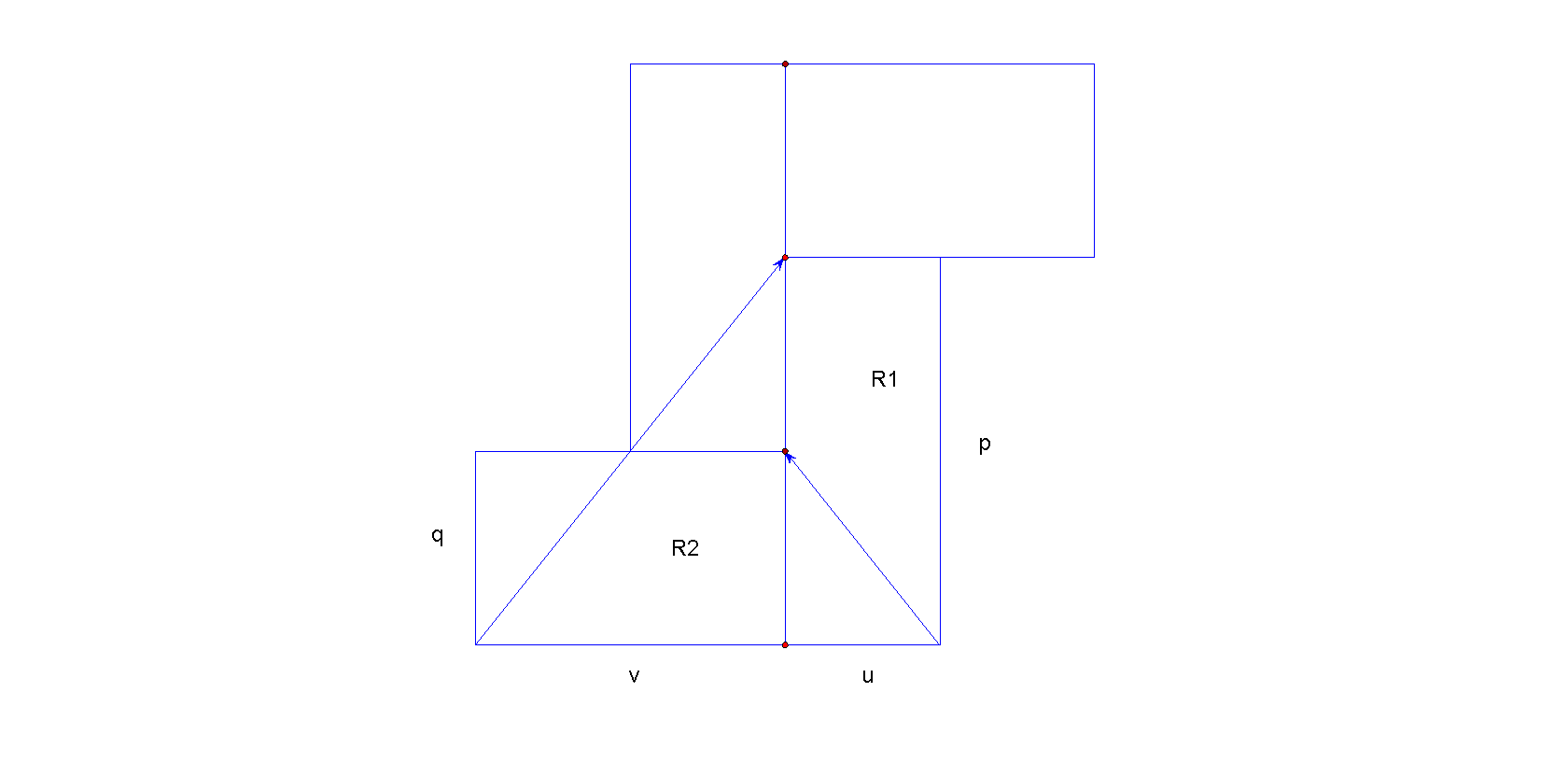}
\end{figure}
with horizontal and vertical sides.

Let us consider two rectangles  $R_1, R_2 \subset \mathbb R^2$
in the position shown in Fig.~1.
 Let the lengths of their horizontal
sides be $u,v$, and of their vertical sides $p,q$.

We consider the lattice of translations $L$ generated by the two
vectors $[v,p]$ and  $[-u,q]$, Fig.~1.

Dividing the plane by the action of the group of translations
gives us the torus $
\mathbb T^2 = \mathbb R^2/L
$
and the natural projection
$
\pi : \mathbb R^2 \to \mathbb R^2/L = \mathbb T^2
$

Translations of $R_1$ and $R_2$ by the vectors from $L$ tile the plane.
The union $R_1 \cup R_2$ is the fundamental domain of the torus.

\begin{definition}\rm\label{defroz}
  The partition $\{R_1,R_2\}$ is called a {\it bi-partition}
of the torus.
A rectangle $R_1$ of a bi-partition $\{R_1,R_2\}$
is called
an {\it isolated rectangle} if it projects $1-1$
under the projection
$\pi : R^2 \to \mathbb R^2/L = \mathbb T^2$.
Otherwise it is called
a {\it connected rectangle}.
\end{definition}

At least one of the rectangles in a bi-partion is connected.
The $(u,p)$ rectangle is isolated iff $u < v, p < q$.
A bi-partiton is called {\it connected} if both of its rectangles are connected,
otherwise it is called {\it isolated}. A bi-partition is isolated iff
$(u -v)(p-q) \leq 0$.

\begin{definition}\rm\label{szkielet}
For a given bi-partition the union of the horizontal sides, $J^s \subset
\mathbb T^2$,
of the rectangles is called the
{\it horizontal spine}, the union of the vertical sides, $J^u \subset
\mathbb T^2$,
is called the {\it vertical spine}.
\end{definition}
The horizontal and the vertical spines are intervals in
$\mathbb T^2$ which intersect in four points
$J^s\cap J^u$, Fig.1.

Let us now reverse the question, given two transversal irrational directions
in the torus, let us find all bi-partitions with the sides of the rectangles
having these directions.
In other words, let $L \subset \Bbb R^2$ be a lattice of translations
isomorphic to $\Bbb Z^2$, and with neither horizontal nor vertical translations,
i.e., the lattice $L$ has no nonzero elements in the coordinate axes.
We are looking for two rectangles
$R_1, R_2,$ which form the bi-partition of the torus $\Bbb R^2/L$.
It is clear from the above construction that such bi-partitions are,
one-to-one, associated with bases $\{ e, f\}$ of the lattice $L$,
such that $e$ belongs to the first, and $f$ belongs to the second quadrant.
Indeed, if $e = [v,p]$ and  $f = [-u,q]$,  then the rectangles $R_1, R_2,$
with the horizontal sides equal to, respectively, $u$ and $v$,
and vertical sides equal to $p$ and $q$, Fig. 1, give us the bi-partition.

Let us consider the family $\mathcal F$ of such bases of $L$.
$\mathcal F$ is always nonempty.
Indeed, let $\{a,b\}$ be a basis in $L$. One of the four bases
$\{\pm a,\pm b\}$ has the property that
the first element, with respect to
the positive orientation, is in the right, and the second element in the left
half-plane.   Let us denote  such a basis by $\{ a_0,b_0\}$.
Now we construct inductively a sequence of such bases of $L$
by the following {\it cutting algorithm}, which is reminiscent of the
coding from the paper of Series, \cite{Se}. Given the basis
$\{ a_n,b_n\}$, we consider $c_n = a_n+b_n$. If $c_n$ is in the right
half-plane then $a_{n+1} = c_n, b_{n+1} = b_n$, and if
$c_n$ is in the left
half-plane then $a_{n+1} = a_n, b_{n+1} = c_n$. The cutting algorithm must
deliver a basis in $\mathcal F$,
since the directions of $a_n$, and $b_n$, converge to
the vertical.

Given a basis $\{e_0,f_0\} \in \mathcal F$ we construct the sequence
of bases $\{e_n,f_n\} \in \mathcal F, n = 0,1,2,\dots$
by the cutting algorithm. Moreover we can run the algorithm
backwards, i.e., for a basis $\{e_n,f_n\} \in \mathcal F$, we consider
the elements of $L$: $e_n-f_n$ and  $f_n-e_n$. One, and only one of them
is in the upper half-plane, let us denote it by $g_n$. If $g_n$ is
in the second quadrant then the basis $\{e_{n-1}, f_{n-1}\}\in \mathcal F$,
where $e_{n-1} = e_n, f_{n-1} = g_n$. If $g_n$ is
in the first quadrant then we put $e_{n-1} = g_n, f_{n-1} = f_n$,
and still get a basis in $\mathcal F$.

Let us summarize our construction: given a basis
$\{e_0,f_0\} \in \mathcal F$ we obtained a series of bases
$\{e_n,f_n\} \in \mathcal F, n = 0,\pm 1, \pm2, \dots$.
We claim that
$\mathcal F$ contains exactly the sequence, and nothing else.
To prove that let us introduce an ordering in $\mathcal F$.
\begin{definition}\rm\label{porzbaz}
A basis  $\{\hat{e},\hat{f}\}\in \mathcal F$ \textit{succeeds} another basis
$\{e,f\} \in \mathcal F$, which we denote by $\{\hat{e},\hat{f}\} \succ \{e,f\}  $,
if $\hat{e},\hat{f}$ are linear combinations of $e,f,$ with nonnegative
coefficients.
\end{definition}

The ordering of $\mathcal F$ is linear, i.e., for any two
bases $\{e,f\}, \{\hat{e},\hat{f}\} \in \mathcal F$,
either $\{e,f\} \succ  \{\hat{e},\hat{f}\}$ or
$ \{\hat{e},\hat{f}\}\succ \{e,f\}$. Indeed let us recall that
our lattice $L$ is isomorphic with $\mathbb Z^2$. Without loss of generality
we can consider $\{e,f\}$ to be the standard basis in $\mathbb Z^2$,
and the horizontal axis becomes a line with the negative slope. The other
basis $ \{\hat{e},\hat{f}\}$ is on the same side of this line as
the standard basis $\{e,f\}$.
\begin{lemma}\label{polbaz}
Let $\{e,f\}$ be the standard basis in $\mathbb Z^2$, and let
$ \{\hat{e},\hat{f}\}$ be another basis. If $\hat{e}$ has both
positive coordinates (i.e., it lies strictly in the first
quadrant), then $\hat{f}$ lies either in the first or third quadrant.
\end{lemma}
It follows from this Lemma that for any two
bases $\{e,f\}, \{\hat{e},\hat{f}\} \in \mathcal F$, if $\hat{e}$
is a linear combination of $e$ and $f$ with positive coefficients then
$ \{\hat{e},\hat{f}\}\succ \{e,f\}$. It is the main step in proving that
the ordering of $\mathcal F$ is isomorphic to
the ordering of $\mathbb Z$. In the following we will call
$\mathcal F$ {\it a fan of bi-partitions, or a bi-fan}.

With a fixed basis $\{e_0,f_0\} \in \mathcal F$
the bi-fan is completely described by {\it the cutting sequence}
$\{s_n\}$ defined as follows: for $n \in \mathbb Z$
\[
s_n = 0 \ \text{if} \ e_{n+1}= e_{n},  \ \ s_n = 1 \ \text{if} \
f_{n+1}= f_{n}.
\]
The bi-fan $\mathcal F$
can be recovered from the basis $\{e_0,f_0\}$ and the cutting
sequence by the following formula: for any $k < n$
\begin{equation}
\label{eq: Z}
[e_{n},f_{n}] = [e_k,f_k]
\left[ \begin{array}{rr} 1 & 0 \\ 1 & 1  \end{array} \right]^{s_{k}}
 \left[ \begin{array}{rr} 1 & 1 \\ 0 & 1  \end{array} \right]^{1-s_{k}}
\dots
\left[ \begin{array}{rr} 1 & 0 \\ 1 & 1  \end{array} \right]^{s_{n-1}}
 \left[ \begin{array}{rr} 1 & 1 \\ 0 & 1  \end{array} \right]^{1-s_{n-1}}
\end{equation}
where $[e_n,f_n]$ is understood as the matrix with columns equal to
$e_n,f_n,$ respectively.

\begin{definition}\rm\label{bazazred}
A basis  $\{e_n,f_n\}\in \mathcal F$ is called a {\it reduced} basis
of the bi-fan if $s_{n-1} \neq s_n$. Otherwise  it is called
{\it intermediate}.
\end{definition}
The terms reduced and itermediate are motivated by
the reduced continued fractions,
and intermiediate continued fractions, due to the connections
of the cutting sequence and continued fractions (cf. \cite{S},\cite{S-W}).
\begin{proposition}\label{bazazredroz}
A basis  $\{e_n,f_n\}\in \mathcal F$ is reduced if and only if the corresponding
bi-partition is isolated.
\end{proposition}
\begin{Proof}
Let $e_n = [v,p]$ and  $f_n = [-u,q]$. We have $s_n = 0$ if and only if
$e_n + f_n$ is in the {\bf left} half-plane, i.e., $v-u < 0$. Further $s_{n-1} = 1$
if and only if the vector $e_n-f_n$ is in the {\bf upper} half-plane, i.e.,
$p-q > 0$. This gives us the condition $(u-v)(p-q) > 0$.
The case $s_n = 1, s_{n-1}=0$ is  characterized by the same inequality.
\end{Proof}$\square$

\begin{figure}[h]{Fig. 2}
  \centering
    \includegraphics[width=0.9\textwidth]{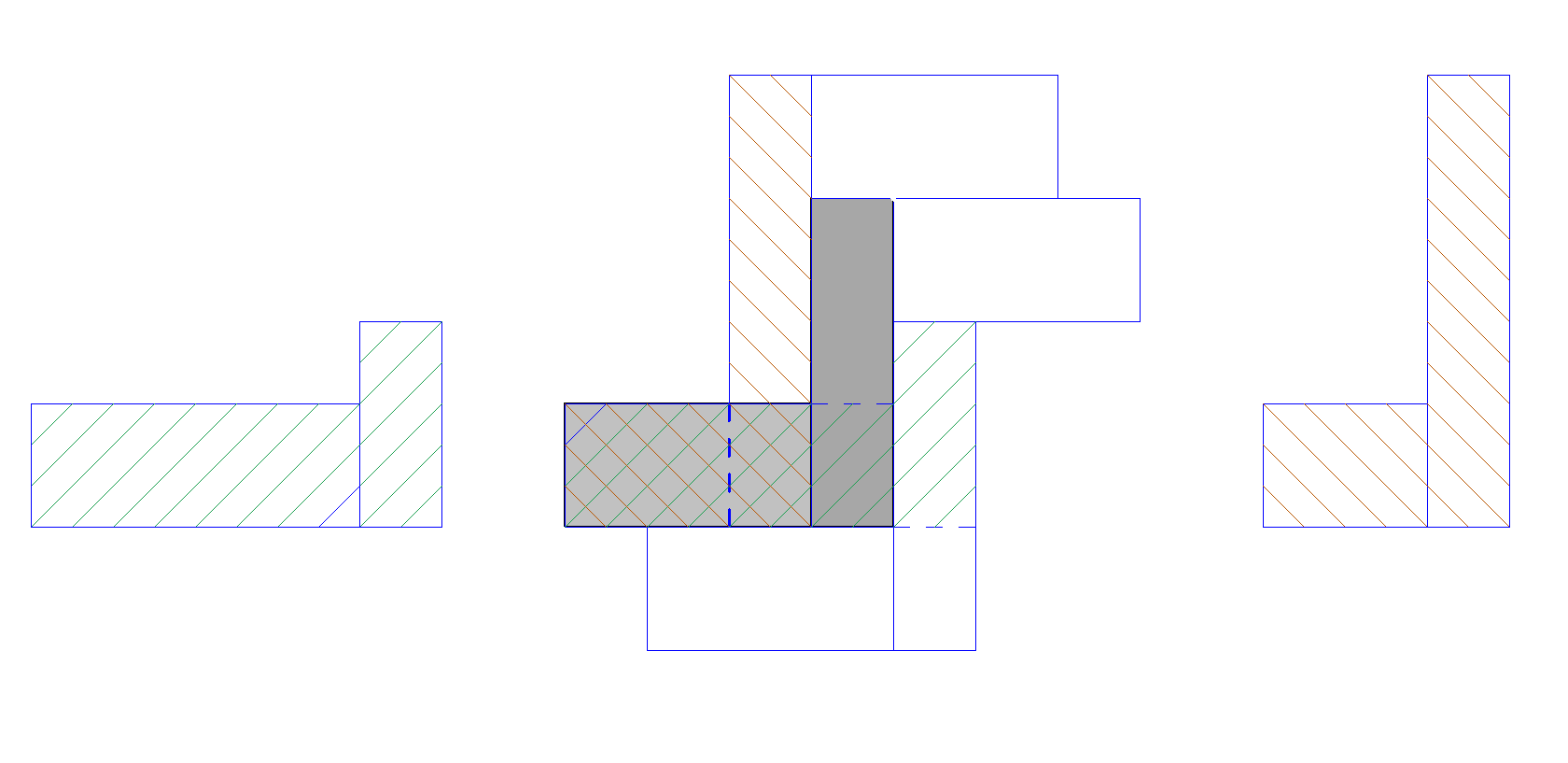}
\end{figure}
The ordering of the bases in the bi-fan gives rise to the respective
ordering of bi-partitions (or rather their shapes), which was used
by Snavely, \cite{Sn}. This ordering
is geometrically transparent, in Fig.~2 a bi-partition is shown with
its successor and predecessor.

\section{Hyperbolic toral automorphism and its bi-fan}\label{hiptoraut}

Let us  consider a hyperbolic toral automorphism $\mathcal D$.
It has a stable and unstable directions.
We choose them as the horizontal axis and the vertical axis, respectively.
With such a choice of coordinates
our automorphism is described by a diagonal matrix $D$
\[
D = \left[ \begin{array}{rr} \mu & 0 \\
0 & \lambda  \end{array} \right]
\ \ \ |\lambda| > 1, |\mu| = \frac{1}{|\lambda|} < 1
\]
We denote the lattice of deck translations
by $L \subset \mathbb R^2$.

Stable and unstable lines have irrational directions. In Section~\ref{bipart} we
associated with such a pair of lines the bi-fan $\F$ of bases in $L$
with vectors lying in the first and the second quadrant, respectively.
Each such basis gives rise to a bi-partition of the torus
$\mathbb T^2 = \mathbb R^2/L$.

\begin{proposition}\label{central}
For any automorphism $\mathcal B$ of the torus which preserves the vertical
and the horizontal axes, we have that $\pm\mathcal B(\mathcal F) = \mathcal F$,
where the sign is the same as the sign of the eigenvalue of  $\mathcal B$
associated with the vertical eigenvector.
Moreover there is an integer $K$ such that for any $\{e_n,f_n\} \in \mathcal F,
\pm\mathcal B( \{e_n,f_n\}) = \{e_{n+K},f_{n+K}\}$,
and for the cutting sequence we get $s_{n+K} = s_{n}$
if $\mathcal B$ is orientation preserving, and  $s_{n+K} = 1-s_{n}$ if
it is orientation reversing.
\end{proposition}
\begin{Proof}
The automorphism $\pm\mathcal B$ maps the vertical line into itself,
preserving its orientation. Hence it either maps
the first and the second quadrants into themselves if $\B$ is
orientation preserving, or exchanges them if it reverses the orientation.
In either case a basis from the
bi-fan $\mathcal F$ is mapped into another basis from $\mathcal F$.
The mapping is 1-1, and moreover the order in $\mathcal F$ has to be preserved.
It follows that there must be an integer $K$ such that
$\pm\mathcal B(e_n) = e_{n+K}, \pm\mathcal B(f_n) = f_{n+K}$
if $\mathcal B$ is orientation preserving, and
$\pm\mathcal B(e_n) = f_{n+K}, \pm\mathcal B(f_n) = e_{n+K}$
if it is orientation reversing. The last property follows immediately.

Note that unless $\B =\pm I$ it must be a hyperbolic toral automorphism.
Further
the integer $K$ is negative if the vertical line is stable
for $\B$.
\end{Proof}$\square$

The above proposition shows that the cutting sequence of our pair
of lines must be periodic. We need a more general definition.
\begin{definition}\rm\label{okres}
A natural number $N$ is called a (\textit{semi-)period} of the cutting sequence
if for every $k$
\[
(s_{k+N} = 1-s_k) \ \ s_{k+N} = s_k.
\]
The sequence is called (semi-)periodic  if it has a (semi-)period.
\end{definition}

Clearly if $N$ is a semi-period of a sequence then $2N$ is a period.
However periodic sequences in general are not semi-periodic.
As usual we consider the smallest (semi-)period and call it
the basic (semi-)period.

It turns out that the (semi-)periodicity of the cutting sequence characterizes
the stable and unstable lines of hyperbolic automorphisms, among all pairs
of irrational directions. It is a consequence of the
following converse of Proposition~\ref{central}.

Let us consider the basic semi-period $N$, or the period if a semi-period is not
present, of the cutting sequence.
\begin{theorem}\label{centralopis}
The toral automorphism $\G$ defined by
$\G(e_0) = f_{N}, \G (f_0) = e_{N}$
in the case of the semi-period, and
$\G (e_0) = e_{N}, \G(f_0) = f_{N}$
in the case of the period,
preserves the vertical and the horizontal axes,
and there is a natural $K$ such that $\G^K = \pm\mathcal D$.

Moreover the centralizer
of the hyperbolic toral automorphisim $\mathcal D$
in $GL(2, \mathbb Z)$ is equal to $\{\pm \G^k| k\in \mathbb Z\}$.
\end{theorem}
\begin{Proof}
The first observation is the following.
Let $\{a,b\}$ be a basis in the lattice $L$. For any irrational
line between $a$ and $b$ we define the infinite one sided forward
cutting sequence by the cutting algorithm.
Two different irrational lines must have different sequences.

By the (semi-)periodicity of the cutting sequence
the forward cutting sequences of the vertical line with respect to the bases
$\{e_0,f_0\}$ and  $\{e_{N},f_{N}\}$ are the same. It follows that
$\G$ takes the vertical line into vertical line. Further,
any basis $\{e_n,f_n\}\in \mathcal F, n\geq 0$ is mapped by $\G$ into
another basis from the bi-fan, with the preservation of the order in
$\mathcal F$.
Hence $\G(\{e_n,f_n\}) = \{e_{n+N},f_{n+N}\}$ for $n\geq 0$, and
for the natural $K$ such that
$\pm\D(\{e_0,f_0\}) = \{e_{KN},f_{KN}\}$ we
get $\G^K = \pm\mathcal D$. By necessity the mapping $\G$
is a hyperbolic toral automorphism with the horizontal stable direction,
and vertical unstable direction.

Any toral automorphism   $\A$ commuting with
$\mathcal D$ preserves the horizontal and the vertical directions,
and hence by Proposition~\ref{central} it leads to a (semi-)period $M$ of the
cutting sequence. This (semi-)period must be a multiple of the basic
(semi-)period $N$, which leads to $\A = \pm\G^{s}$
for an integer $s$ such that $M = |s|N$. The sign of $s$ depends on the
automorphism $\A$ having vertical stable, or unstable direction.
\end{Proof}$\square$

In the following we will refer to the mapping $\G$ as {\it the generator}
of $\D$.

It follows from \eqref{eq: Z} that in the basis $\{e_k,f_k\}$
the generator $\G$ is represented by the following
matrix with nonnegative entries
\begin{equation}\label{factorization}
\left[ \begin{array}{rr} 1 & 0 \\ 1 & 1  \end{array} \right]^{s_{k}}
 \left[ \begin{array}{rr} 1 & 1 \\ 0 & 1  \end{array} \right]^{1-s_{k}}
\dots
\left[ \begin{array}{rr} 1 & 0 \\ 1 & 1  \end{array} \right]^{s_{k+N-1}}
 \left[ \begin{array}{rr} 1 & 1 \\ 0 & 1  \end{array} \right]^{1-s_{k+N-1}}
\left[ \begin{array}{rr} 0 & 1 \\ 1 & 0  \end{array} \right]^a
\end{equation}
where $a =1$ if $N$ is the basic semi-period, and $a=0$ otherwise.
If we change the order of elements in the basis the representation
will change by the conjugation by
$\sigma =\left[ \begin{array}{rr} 0 & 1 \\ 1 & 0  \end{array} \right]$.
Hence we get $2N$ different matrices with nonnegative entries conjugate to
our toral automorphism $\D$. It turns out that there are no other
such matrices.

\begin{proposition}\label{wachlarz}
If the automorphism $\D$ (or its generator $\G$)
is represented in a basis
$\{a,b\}$ of $L$ by a matrix with nonnegative
elements then either $\{a,b\}\in \F$ or  $\{-a,-b\}\in \F$.
\end{proposition}
\begin{Proof}
A hyperbolic toral automorphism defined by a matrix with nonnegative elements
has the unstable direction with positive slope, and the stable
direction with negative slope. This means that the vectors of the
standard basis lie on different sides of the unstable line, and
on the same side of the stable line, which proves our claim.
\end{Proof}$\square$

The factorization (\ref{factorization}) appears explicitely in the paper of
Appelgate and Onishi, \cite{A-O},
see also \cite{S-W}.

\section{Berg partitions of hyperbolic toral automorphisms}\label{Bergpartdlaaut}

We now turn to Markov partitions of $\mathcal D$ with two rectangular
elements, which are hence bi-partitions.
\begin{definition}\rm\label{rozbicieBerga}
A bi-partition $\{R_1,R_2\}$ with the spines $J^s$ and $J^u$ is a \textit{Berg
partition} of the hyperbolic toral automorphism $\mathcal D
:\mathbb R^2/L \to \mathbb R^2/L$
if
\[
\mathcal D (J^s) \subset J^s, \ \ \ \ \mathcal D (J^u) \supset J^u.
\]
\end{definition}
It follows from this definition that the spines have stable and unstable
directions, and contain
fixed points of $\mathcal D$.

The discussion in Section~\ref{bipart} leads to the conclusion that for each Berg
partition $\{R_1,R_2\}$ there is an element $\{e_0,f_0\}$
of the bi-fan $\mathcal F$ associated with it as in Fig 1.

Every Berg partition comes with the connectivity matrix
\[
C = \left[ \begin{array}{rr} k & l \\ m &  n \end{array} \right].
\ \ \ k,l,m,n \ \in \ \mathbb Z^+.
\]

The element $c_{ij}$
in $i$-th row and $j$-th column of the connectivity matrix $C$
is the number of translates of $R_j$ in the
plane $\mathbb R^2$, intersected by the image $\D R_i$, lifted
to the plane.

Adler \cite{Ad}, and  Manning \cite{M} , showed that if
a hyperbolic toral automorphism is conjugate to $\pm C^T$, for a matrix
$C$ with nonnegative entries, then it has
a 2-element Markov partition
with the connectivity matrix $C$. One of the byproducts of our approach
is a simple proof of this fact.

\begin{theorem}\label{repaut}
For a Berg partititon \{$R_1,R_2$\} of $\D$ with the connectivity matrix  $C$,
the automorphism $\D$ is represented in the basis $e_0= [v,p],f_0=[-u,q]$
by the matrix $\pm C^T$, where the sign is equal to the sign of the
trace of $\D$ (i.e., it is plus if $\lambda > 1$,
and minus, if $\lambda < -1$).
\end{theorem}
\begin{Proof}
The eigenvalues of $\D$ are $\lambda$ and $\mu$.
The images of the rectangles $R_1,R_2$ are stretched vertically by the factor
$|\lambda| > 1$, and contracted horizontally by the factor $|\mu| < 1$.
These images intersect completely (i.e., from top to bottom)
certain rectangles in the tiling. This gives us the following
``covering conditions''
\[
\begin{aligned}
&|\lambda|p = k p + l q \\
&|\lambda|q = m p + n q
\end{aligned},
\]
for some nonnegative integers $k,l,m,n$.
Further the translates of these images
fill exactly both $R_1$ and $R_2$ (or any element of the tiling). Hence we get
the following ``packing conditions'':
\[
\begin{aligned}
&u = k|\mu|u + m|\mu| v \\
&v = l|\mu|u + n|\mu| v
\end{aligned}
\]

We conclude that $|\lambda|$ is the Perron-Frobenius
eigenvalue of $C$, and the vector $(p,q)$ is the respective column
eigenvector.
Since $|\mu| = 1/|\lambda|$, the vector $(u,v)$ is the
row eigenvector of $C$ with the same eigenvalue.

Let $\nu$ denote the other
eigenvalue of $C$. It is straightforward that the vector $(v,-u)$ is the
column eigenvector of $C$ with eigenvalue $\nu$ (and the vector $(q,-p)$ is
the respective row eigenvector). We can put this together into
\begin{equation}
\label{eq: C}
\left[ \begin{array}{rr} \nu & 0 \\ 0 & |\lambda|  \end{array} \right]
\left[ \begin{array}{rr} v & -u \\ p &  q \end{array} \right]
=
\left[ \begin{array}{rr} v & -u \\ p &  q \end{array} \right]C^T.
\end{equation}

The automorphism $\D$ is represented in the basis $\{e_0,f_0\}$
by an integer matrix $F \in GL(2,\mathbb Z)$ which translates into
\begin{equation}
\left[ \begin{array}{rr} \mu & 0 \\ 0 & \lambda  \end{array} \right]
\left[ \begin{array}{rr} v & -u \\ p &  q \end{array} \right]
=
\left[ \begin{array}{rr} v & -u \\ p &  q \end{array} \right]F.
\end{equation}

Comparing (3) and (4) we conclude that the $2\times 2$ integer
matrices $F$ and $C^T$ have the same row eigenvectors
(equal to $(v,-u)$ and $(p,q)$). Hence also the integer matrix $F^{-1}C^T$
has the same eigenvectors, and its eigenvalues are
$|\lambda|/\lambda = \pm 1$ and $\nu/\mu$. Should these
eigenvalues be different, the eigenvectors of $F^{-1}C^T$ would have
rational slopes. Hence $\nu/\mu = |\lambda|/\lambda = \pm 1$ and
we obtain $F = \pm C^T$.
\end{Proof}$\square$

In general a bi-partition is not a Berg partition.
However every bi-partition with sides parallel to stable
and unstable directions of $\D$
is a translate of a Berg partition. Indeed a bi-partition can be
translated so that one of the four intersection points of the horizontal
and vertical spines is a fixed point. Such a bi-partiton is then by
necessity a Berg partition in the case of both positive eigenvalues.

More generally for  a bi-partition,
in the case of both positive eigenvalues $\lambda,\mu$,
if the horizontal and vertical spines contain fixed points then
we have a Berg partition. In the case of a negative eigenvalue
it is not so, because if the fixed point is too close to the endpoint of
the respective spine, then the condition from Definition~\ref{rozbicieBerga} is not
satisfied. This leads us to
\begin{definition}\rm\label{middle}
For any $\beta >1$
the $\beta$-{\it middle} $c(J)$ of a segment $J$ is the middle subsegment
with the length
$
|c(J)| = \frac{\beta -1}{\beta +1}|J|
$
\end{definition}
It can be checked by direct calculation
(and it will be done in Section~\ref{Bergpartdlamacierzy})
that all the four common points
of the horizontal and vertical spines of a bi-partition
of a hyperbolic toral automorphism with eigenvalues $\lambda,\nu$
are either endpoints of a spine, or lie in the $|\lambda|$-middle of the spine.

Let us now consider a hyperbolic toral automorphism $\D$,
and a bi-partition with the stable and unstable directions
of the spines. In view of the preceding discussion we get
\begin{theorem}\label{rozbdlaD}
If the spines $J^s$ and $J^u$,
or their
$|\lambda|$-middles $c(J^s)$ and $c(J^u)$,
contain fixed points of $\D$,
depending on the sign of the eigenvalues $\lambda, \mu$,
then the bi-partition is a Berg partition for $\D$.
\end{theorem}
Except for the case of both negative eigenvalues,
it follows from Theorem~\ref{rozbdlaD} that every hyperbolic automorphism
has Berg partitions.
Indeed in every other case we can place a fixed point of $\D$ in
one of the four intersection points of the horizontal and vertical spines,
and not at the endpoint of the spine, if the respective eigenvalue is negative.
Since the common point is guaranteed to be in the $|\lambda|$-middle of the
respective spine, the conditions of Theorem~\ref{rozbdlaD} are satisfied.

If both  eigenvalues are negative, we need to use two different fixed points.
It turns out that every two fixed points can be used.
In the next section we will count the number of nonequivalent Berg
partitions with a given connectivity matrix. It will deliver the
existence of such partitions also in the negative case.

Let us note that the theorem of Adler and Manning will follow from our
discussion. Indeed if the automorphism
$\D$ is represented in a certain basis
by $\pm C^T$, for a matrix $C$ with nonnegative
entries, then by Proposition~\ref{wachlarz} the basis (or its equivalent)
must belong to the bi-fan. The corresponding bi-partition can be then
translated into a Berg partition. By Theorem~\ref{repaut}
the connectivity matrix of this partition must be equal to $C$.

\section{Berg partitions with a fixed connectivity matrix $C$}\label{Bergpartdlamacierzy}

We are going to count the number of different Berg partitions
with the fixed connectivity matrix $C$,
for a fixed toral automorphism $\D$, conjugate to $C^T$ or
$-C^T$.

\begin{definition}\rm\label{centgrupa}
 The {\it centralizer} $Z(\D)$
is the group of homeomorphisms of the torus commuting with
$\D$.
\end{definition}
It was proven by Arov \cite{Ar}, and Adler and Palais
\cite{A-P} that the centralizer $Z(\D)$ contains only
affine maps. By Theorem~\ref{centralopis} all the toral automorphisms
in $Z(\D)$ have the form $\pm \G^k$ for the generator $\G$
and some integer $k$. The centralizer $Z(\D)$ is generated
by these, and the translations which take zero to other fixed
points of $\D$.

\begin{definition}\rm\label{rozbrow}
Two Berg partitions are {\it equivalent} if there is
a mapping in $Z(\D)$ which takes one
partition into the other,
without regard to the ordering of the rectangles.
\end{definition}

We want to find the number of equivalence classes of
Berg partitions for a given toral automorphism $\D$.

To each Berg partition we associated in Section~\ref{Bergpartdlaaut}
a basis in the respective bi-fan $\F$.
It is enough to consider Berg partitions associated
with $N$ consecutive bases of $\F$,
$\{e_1,f_1\},\dots,\{e_N,f_N\}$, where $N$ is the (semi-)period
of the cutting sequence. By Theorem~\ref{repaut}
a basis in $\F$ determines the connectivity matrix of the
respective Berg partition. In particular, as it was proven
in Section~\ref{hiptoraut}, the connectivity
matrices for the $N$ consecutive bases are all different.

Let us note that for a given basis in $\F$ we get two
different connectivity matrices $C$ and $\sigma C \sigma$
where $\sigma =\left[ \begin{array}{rr} 0 & 1 \\ 1 & 0  \end{array} \right]$.
However this does not lead to nonequivalent Berg partitions
since such a change in the connectivity matrix is effected
by merely reordering the elements of the Berg partition.

Berg partitions that have connectivity matrices differing
by more than the conjugation by $\sigma$ are nonequivalent.

Hence we get $N$ nonequivalent shapes of Berg partitions.
Moreover every Berg partition is equivalent to one of those,
which is the content of the following

\begin{theorem}\label{przesrozb}
If $\{R_1,R_2\}$ and $\{R'_1,R'_2\}$ are two Berg partitions with the
same connectivity matrix then there is a translation $\mathcal T$ of
$\T^2$
such that $\{\mathcal T(R_1),$ $\mathcal T(R_2)\}$
is equivalent to $\{R'_1,R'_2\}$.
\end{theorem}
\begin{Proof}
Without loss of generality we can assume that
$\{e_0,f_0\} \in \F$ is the basis associated with the
bi-partition $\{R_1,R_2\}$, and $\{e_m,f_m\} \in \F, m\geq 1,$ is
associated with $\{R'_1,R'_2\}$. Since the two bi-partitions
have the same connectivity matrices, then it follows that
there is a natural $k$ such that $m = kN$. Further
$\G^k$ takes one basis into the other, and hence the bi-partitions
$\{\G^{-k}(R'_1),\G^{-k}(R'_2)\}$ and  $\{R_1,R_2\}$
have the same basis $\{e_0,f_0\} \in \F$ associated to them.
But that means that the rectangles in the bi-partitions are
isometric, and so they differ by a translation.
\end{Proof}$\square$

Let us remark that while translations of a Berg partition may give
us other, nonequivalent, Berg partitions, the reflection
in the vertical axis never delivers one. More precisely, let
$E$ be the reflection in the vertical axis of $\R^2$. If $\{R_1,R_2\}$
is a bi-partition then
$\{E(R_1),E(R_2)\}$ is not a bi-partition of the same torus.
To see this let us consider the basis $\{e,f\} \in \F$
associated with $\{R_1,R_2\}$. Should $\{E(R_1),E(R_2)\}$
be a bi-partition then $\{E(e),E(f)\}$
would also be a basis in $\F$, which is impossible
because it cannot be either earlier or later than
$\{e,f\}$ in the ordering of $\F$.

By Theorem~\ref{przesrozb}, to count the number of nonequivalent
Berg partitions of
$\D$ with the same connectivity matrix we need to choose
a respective bi-partition and
then translate it around the torus. Each time
the horizontal and vertical spines, or their
$|\lambda|$-middles if needed by Theorem~\ref{rozbdlaD},
contain a fixed point we get a Berg partition. And all Berg partitions
are equivalent to one obtained in such a way.
To accomplish our task we need to study the set of fixed points.

Let $\widehat L$ be the superlattice $\widehat L = (D-I)^{-1}L$.
The group $G = \widehat L / L $ is isomorphic to
the subgroup of translations in the centralizer $Z(\D)$.

$G = \widehat L/L$ acts freely and transitively on the set of
fixed points of $\mathcal A$.
In particular the number of fixed points of $\D$
is equal to
\[
|\det (D-I)| =
\begin{cases}
|tr D -2| \ \ \text{if} \ \det D = 1\\
|tr D| \ \ \text{if} \ \det D = -1\\
\end{cases}
\]

Any fixed point can be translated into any other fixed point
by an element from $Z(\D)$. We place one of the fixed points,
$p_1$, into the horizontal spine $J^s$. To get a Berg partition
we now translate the bi-partition by horizontal vectors, keeping
$p_1$ in $J^s$. Each time the vertical spine $J^u$ hits a fixed point
$p_2$ we get a Berg partition, at least if both eigenvalues are positive.
In general, by Theorem~\ref{rozbdlaD} the spines need to be replaced by their
$|\lambda|$-middles, according to the signs of the eigenvalues of $\D$.

At the same time such partitions may be equivalent.
Indeed let us consider the mapping $\mathcal Z$,
the rotation by $\pi$ around $p_1$. Any Berg partition is mapped
by $\mathcal Z$ into an equivalent Berg partition,
which is also a translation by a horizontal vector, see Fig 1.
It is the same bi-partition if
and only if both $p_1$ and $p_2$ are at the centers of the spines.
It happens if and only if the diagonal entries, and the off
diagonal entries of the connectivity matrix $C$ have the same parity.
Indeed, the translation from
the center of the horizontal spine to the center of the vertical
spine is equal to $\frac{1}{2}(e+f)$, and this vector
belongs to the superlattice $\widehat L$ if and only if
$(D-I)(e+f) \in 2L$. Since by Theorem~\ref{repaut} the automorphism
$\D$ is represented in the basis $\{e,f\}$ by $\pm C^T$,
we conclude that the latter is equivalent to the sum of
the columns of $C^T$ having both odd entries. Finally,
using $det~ C = \pm 1$ we get the required claim.

Let the connectivity matrix be equal to
$C=
\left[ \begin{array}{rr} k & l \\ m &  n \end{array} \right]
$. As it was established in Section~\ref{Bergpartdlaaut}, the matrix $C$ is an element in $GL(2,\mathbb Z)$.

\begin{theorem}\label{ilerozb}
There are exactly
\[
\left[\frac{k+l+m+n}{2}\right]
\]
nonequivalent Berg partitions with the connectivity matrix
$C$, where $[\cdot]$ denotes the integer part of a real number.
\end{theorem}
In the proof we will distinguish three cases:

Case 1: $\det C =1, tr \D > 0$,

Case 2: $\det C =-1, tr \D < 0$

Case 3: $\det C = 1, tr \D < 0$.

The reason that we do not need to consider the fourth case,
$\det C = -1, tr \D > 0 $, is the following
\begin{lemma}\label{odwrot}
A Berg partition for a hyperbolic toral automorphism $\D$ with the connectivity
matrix $C$ is also a Berg partition for $\D^{-1}$ with the connectivity matrix
$C^T$.
\end{lemma}
It follows from this Lemma that if $\det C = -1, tr \D > 0$, then the number
of Berg partitions for $\D$ with the connectivity matrix $C$
is equal to the number of Berg partitons for $\D^{-1}$
with the connectivity matrix $C^T$. This gives us the second case:
$\det C^T = -1$ and $tr \D^{-1} < 0$.

\begin{Proof}
Let $\Gamma : J^s\times J^u \to \mathbb R^2$ be the
map defined by $\Gamma(p_1,p_2) = p_2-p_1$.

We need to count the number of ways in which two fixed points can be
placed into the horizontal and vertical spines. This
number is equal to the number of elements in
$\Gamma^{-1}\left(\widehat L\right)$, which is the same as the number of
elements in
$\Gamma(J^s\times J^u)\cap \widehat L$.

In coordinates $\Gamma(J^s\times J^u)$ is the box $P= [-u,v]\times [0,p+q]$.
The generators of $L$ are given by $e'= (D-I)^{-1}e,\;f'= (D-I)^{-1}f,$,
where $e = [v,p], f = [-u,q]$. The position of $e'$, $f'$ with respect to
the box $P$ is different in the three cases.
The problem in counting the number of elements of the lattice $\widehat L$
in $P$ is that the vertices of $P$ do not belong to $\widehat L$.
To get the number  we replace the box $P$ by an appropriate polygon
$Q$
with the vertices in $\widehat L$, and such that the number of elements
of $\widehat L$ in $Q$ is the same as in $P$,
except for the vertices.
In each of the three cases the modified boxes have different geometry.

Once we have a polygon with vertices from our lattice it is
easy to establish the number of elements from $\widehat L$.
For that purpose we observe that
\begin{lemma}\label{rownoleg}
For a closed parellologram in $\mathbb R^2$ with
vertices in  $\mathbb Z^2$, its area is equal to the number
of points from $\mathbb Z^2$ in the interior,
plus one half of their number on the boundary excluding vertices, plus $1$.
\end{lemma}
Further we choose to depict $P$ and $Q$ in coordinates
in which $\widehat L$ becomes $\mathbb Z^2$.
This is easily accomplished based on the formula
\begin{equation}
\label{eq: D}
(A-I)^{-1}
\left[ \begin{array}{rr} v & -u \\ p &  q \end{array} \right](\pm C^T-I)
=
\left[ \begin{array}{rr} v & -u \\ p &  q \end{array} \right].
\end{equation}
Indeed (5) says that if $e'$, $f'$ are the basic vectors then $e,f$ become
the columns of $\pm C^T-I$.
\begin{figure}[h]{Fig. 3.1}
  \centering
    \includegraphics[width=0.80\textwidth]{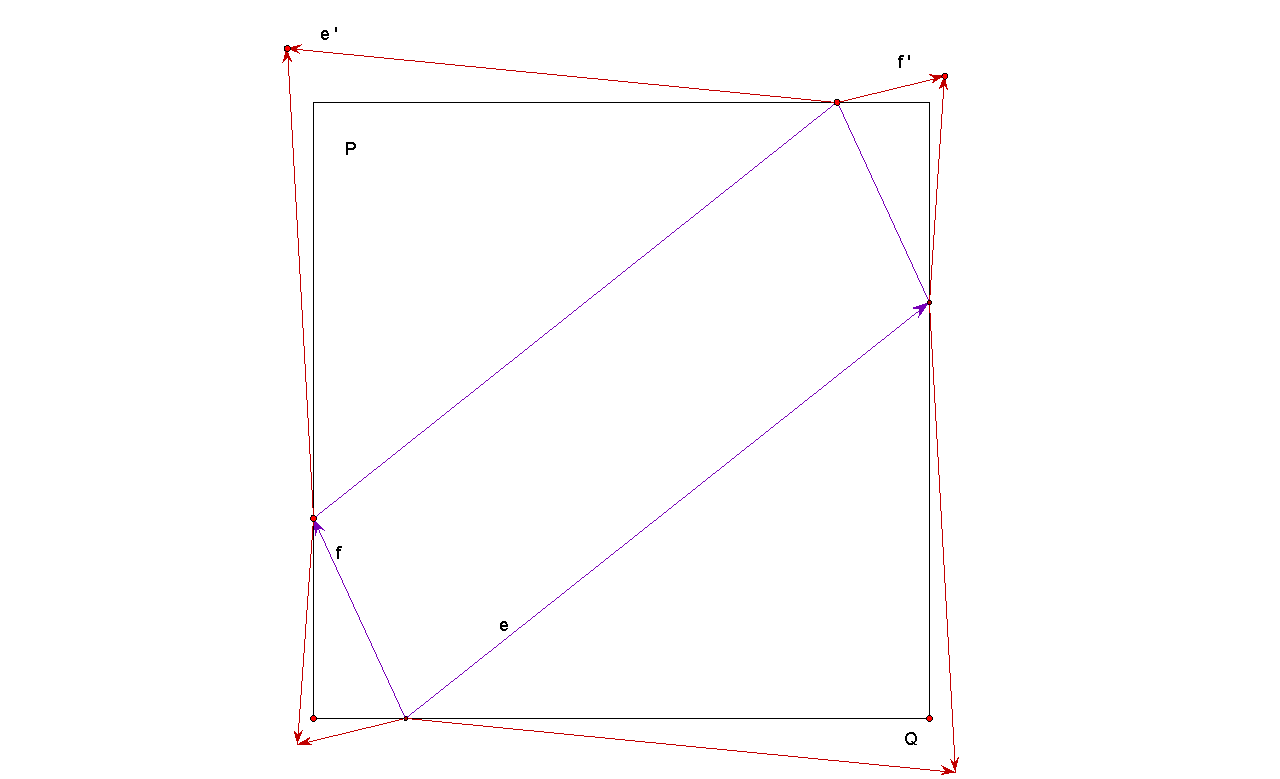}
\end{figure}
Let us denote by $\lambda > 1$ the eigenvalue of the incidence matrix $C$.
The unstable eigenvalue of $\D$ is then equal to $\pm\lambda$
depending on the sign of the trace of $\D$.

In Case 1 we have $e' = [-\frac{\lambda}{\lambda -1}v,\frac{1}{\lambda -1}p],
f' = [\frac{\lambda}{\lambda -1}u,\frac{1}{\lambda -1}q]$
and their position with respect
to the box $P$, and the modified box $Q$ are shown in Fig.~3.1.
In order to establish that $Q$ has no more points from $\widehat L$
than $P$ (except for the four vertices) we pass to the coordinate system
with the basis
$\{e',f'\}$.  The result is depicted in Fig.~4.1.

The crucial feature is that in
view of
\eqref{eq: D} the box $P$ is now
obtained by the following construction.
We take the paralellogram $Y$ spanned by the columns of $C^T$.
There are no integer points in its interior.
The unstable direction is enclosed between the sides of the parallelogram.
We modify the sides to be the columns of $C^T-I$
and get a wider parallelogram which is inscribed
in the box $P$. The box $P$ has the sides with stable and unstable directions.
By inspection we can convince ourselves that the modified box $Q$
has only vertices added to the family of integer points in it, compared to $P$.
Indeed it follows easily from the fact that
the parallelogram $Y$ has no integer points in its interior. The number of integer points in $P$
can be now established to be equal to $k+n+l+m-1$. Each of these integer points gives us a Berg partition.
\begin{figure}[h]{Fig. 4.1}
  \centering
    \includegraphics[width=1\textwidth]{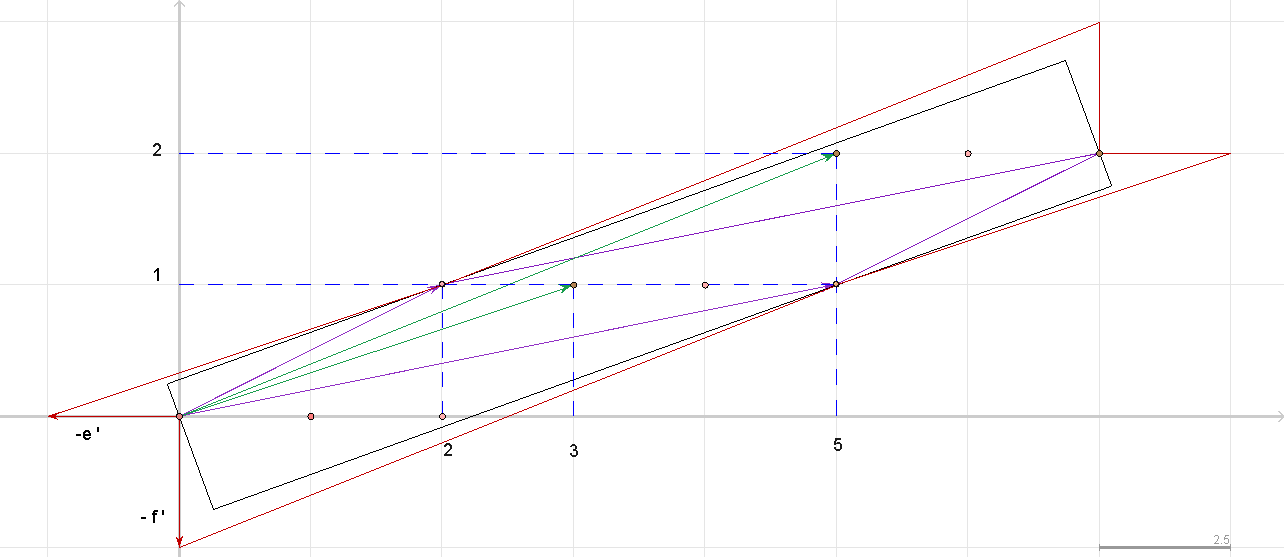}
\end{figure}
And two such Berg partitions are equivalent if and only if the respective integer points are symmetric
with respect to the center of the box $P$.
The center of the box is an integer point if and only if two fixed points can be placed into the centers
of both spines, as it was explained above.

To get the number of nonequivalent Berg partitions we need to divide the number of integer points by $2$,
excepting the center of the box if it is an integer point.
In both cases the result is
$\left[\frac{k+l+m+n}{2}\right]$, where $[\cdot]$ denotes the integer part of a real number.
\begin{figure}[h]{Fig. 3.2}
  \centering
    \includegraphics[width=0.80\textwidth]{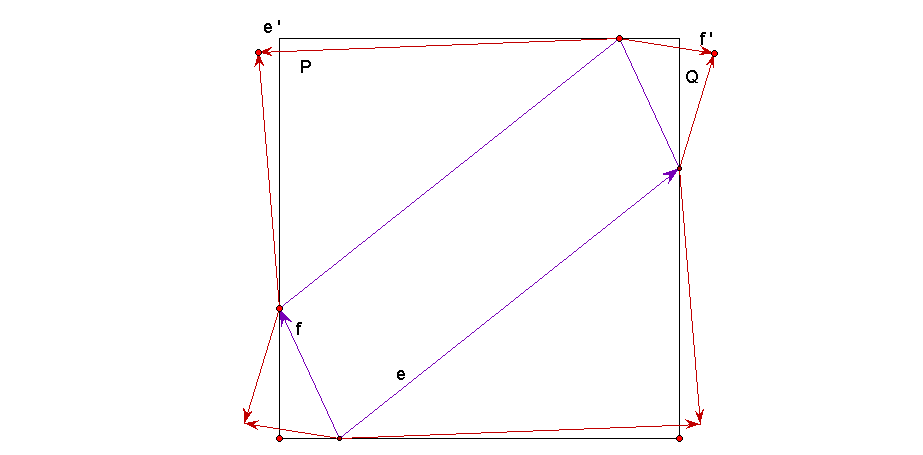}
\end{figure}

In Case 2 we have $e' = [-\frac{\lambda}{\lambda -1}v,-\frac{1}{\lambda +1}p],
f' = [\frac{\lambda}{\lambda -1}u,-\frac{1}{\lambda +1}q]$ and their position with respect
to the box $P$, and the modified box $Q$, are shown in Fig.~3.2.
Further we pass to the coordinate system with the basis
$\{-e',-f'\}$.  The result is depicted in Fig.~4.2.
\begin{figure}[h]{Fig. 4.2}
  \centering
    \includegraphics[width=1\textwidth]{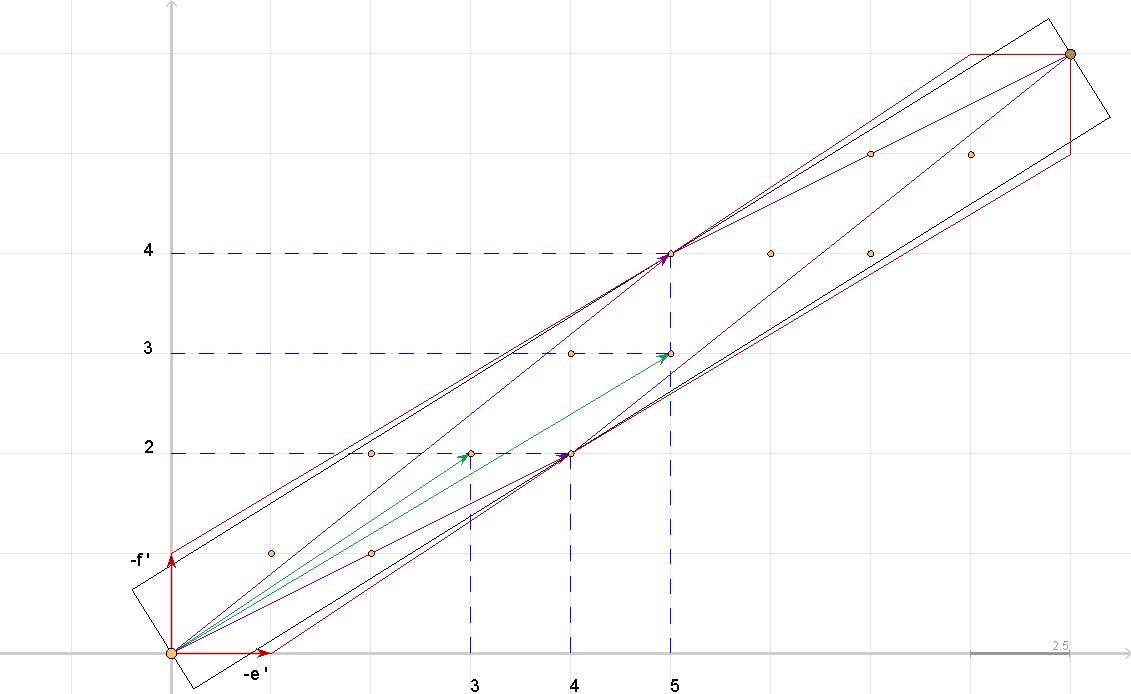}
\end{figure}
 This time we modify the paralellogram $Y$
spanned by the columns of $C^T$ to the parallelogram with the sides equal to the columns of $C^T+I$
and get a wider parallelogram which is inscribed in the box $P$.
Again using the fact that
the parallelogram $Y$ has no integer points in its interior,
we can establish that $Q$ adds only the vertices to the count of integer points.
The result is that the number of integer points in $P$ is equal to $k+n+l+m+1$.
However we need to exclude the points which are too close to the stable boundary of $P$.
It is straighforward that $(\pm(e'+f')$ give us one of the fixed points on the boundary of the $|\lambda|$-middle
of the vertical spine. It follows that we need to exclude only two integer points (on the stable sides of $P$),
and the resulting number of admissible integer points is again  $k+n+l+m-1$.

Finally in Case 3 we have $e' = [-\frac{\lambda}{\lambda +1}v,-\frac{1}{\lambda +1}p],
f' = [\frac{\lambda}{\lambda +1}u,-\frac{1}{\lambda +1}q]$ and their position with respect
\begin{figure}[h]{Fig. 3.3}
  \centering
    \includegraphics[width=0.80\textwidth]{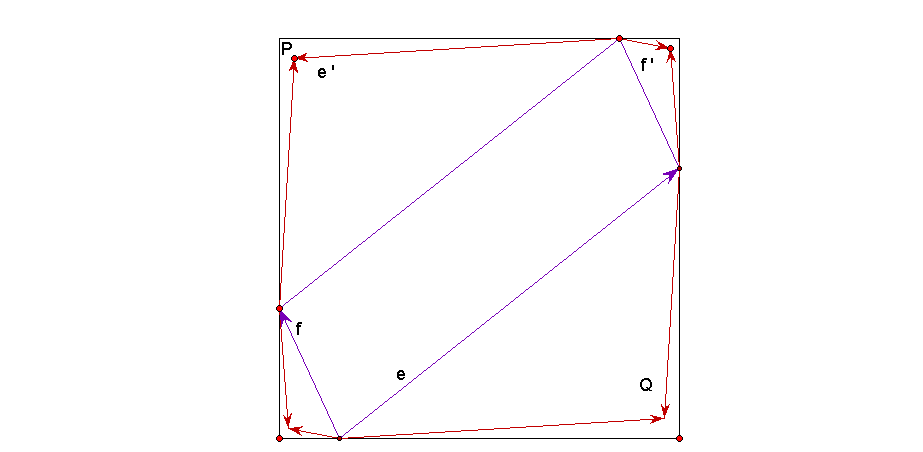}
\end{figure}
to the box $P$, and the modified box $Q$, are shown in Fig.~3.3.
In this case we again pass to the coordinate system with the basis
$\{-e',-f'\}$.  The result is depicted in Fig.~4.3.
\begin{figure}[h]{Fig. 4.3}
  \centering
    \includegraphics[width=1\textwidth]{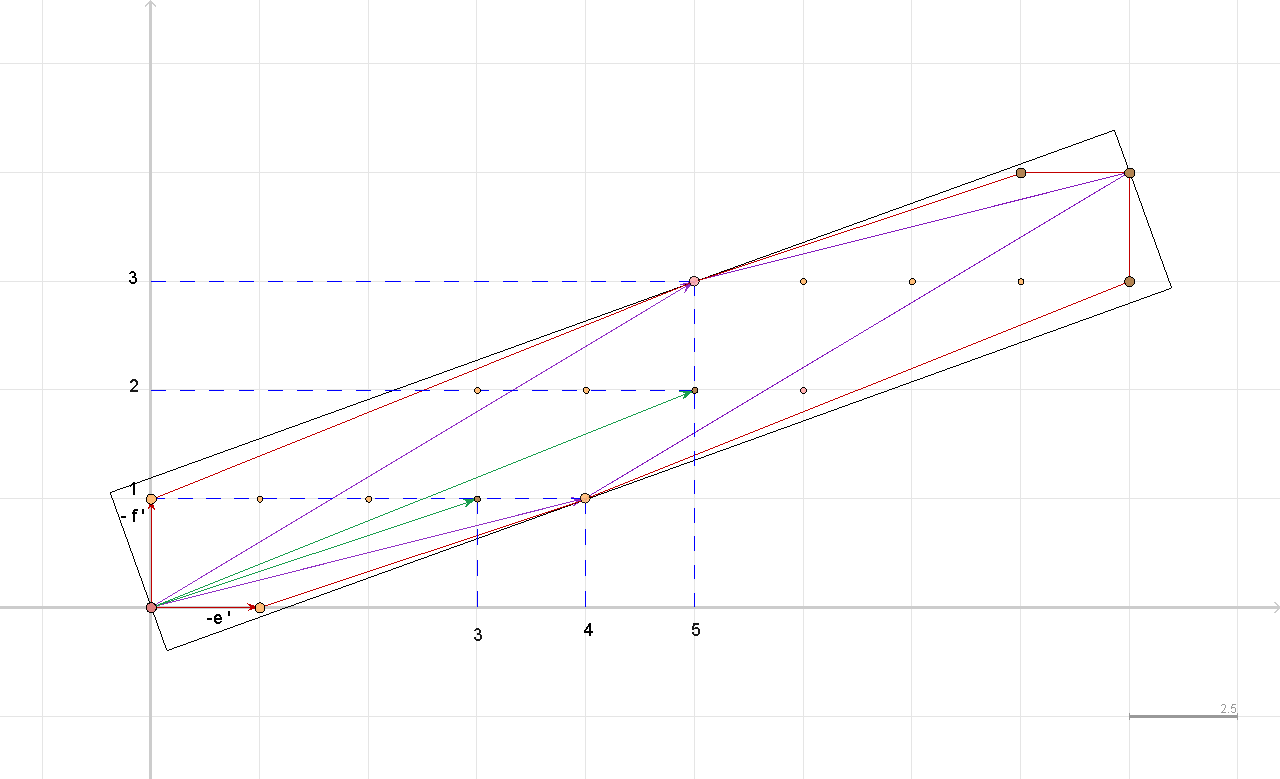}
\end{figure}
 This time the number of integer points in $P$
is equal to $k+n+l+m+7$. But we need to exclude the points which are too close to the  boundary of $P$.
There are four points on the boundary, and the other four points to be excluded are $-e',-f',e+f+e',
e+f+f'$.
The resulting number of admissible integer points is yet again  $k+n+l+m-1$. To prove it we can argue that
by inspection of  Fig.~4.3 there are no more integer points in $P$ which correspond to
a fixed point outside of the $\lambda$-middle of the vertical spine. We still need to check that there are no other
integer points in $Q$ which correspond to a fixed point outside of the $\lambda$-middle of  the horizontal spine.
This can be accomplished by invoking the symmetry used above, that a Berg partition for $\D$ with the connectivity matrix
$C$ is also a Berg partition for $\D^{-1}$ with the connectivity matrix $C^T$. In this symmetry the role of
stable and unstable boundaries of the Berg partition are exchanged.
\end{Proof}$\square$

The above proof of Theorem~\ref{ilerozb} can be developed into a detailed
examination of different Berg partitions with the same connectivity matrix,
and thus with the same shape (cf. Theorem~\ref{przesrozb}).
Two partitions of the same shape (i.e., translates of each other) are
nonequivalent because they have nonequivalent pairs of fixed points
in the respective skeletons, or because the fixed points are placed
differently in the skeletons. The first observation is that any pair of fixed
points of $\D$  can be used to build a Berg partition, with the exception of one
fixed point in the case of both negative eigenvalues (Case 3).
Indeed, every pair of fixed points gets its representation in the paralellogram
$\Sigma \subset P$ spanned by $\{e,f\}$.

We say that a family of Berg partitions differing by translations is
{\it hinged}
if all of them contain the same fixed points in their skeletons,
or that the family is {\it hinged on the fixed points}.
Two translations of the same bi-partition are hinged Berg partitions if
they are represented by two points in the box $P$ differing by a multiple
of vectors $e$, or $f$.

A family of equivalence classes of Berg partitions is said to be
\textit{hinged on a pair of fixed points}
if they contain representatives with these fixed points in their
skeletons.

\begin{proposition}\label{opisilerozb}
A Berg partition with the connectivity matrix
$C=
\left[ \begin{array}{rr} k & l
\\ m &  n \end{array} \right]
$
is isolated if and only if $|m-n| < |k-l|$.

For an isolated bi-partition there can be at most two translates which
are hinged Berg partitions. More precisely among the
$[\frac{k+l+m+n}{2}]$ equivalence classes of isolated Berg partitions
there are exactly $\left[\frac{l}{2}\right] +
\left[\frac{m}{2}\right]$ hinged pairs.

For a connected bi-partition with $l \geq m$ there can be no more than
\[
\left[\sqrt{\frac{l}{m}+\left(\frac{n-k}{2m}\right)^2}-\frac{|n-k|}{2m}
\right] +2
\]
translates which are hinged Berg partitions. This estimate is sharp,
for example for matrices of the form
$C=
\left[ \begin{array}{rr} n-1 & 1 \\ n^2-n-1 &  n \end{array} \right]
$,
for which there are exactly $n+1$ hinged Berg partitions.
\end{proposition}

\begin{Proof}
Expressing  the eigenvectors of $C$ we obtain
\[
\frac{v}{u} = \frac{|\lambda|-k}{m} =\frac{l}{|\lambda|-n}, \ \
\frac{p}{q} = \frac{l}{|\lambda|-k} =\frac{|\lambda|-n}{m},
\]
and the criterion for isolated bi-partitions follows.

\begin{figure}[h]{Fig. 5}
  \centering
    \includegraphics[width=0.8\textwidth]{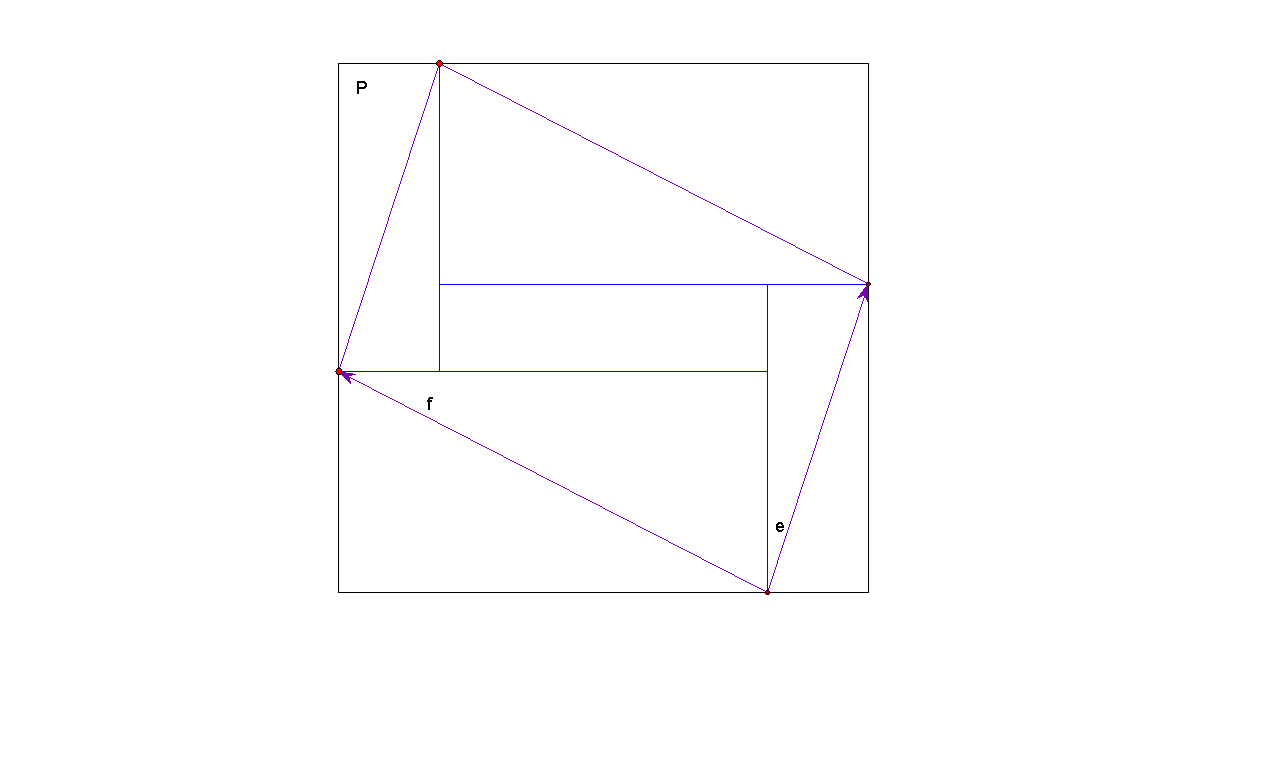}
\end{figure}

To count the number of hinged Berg partitions we need to inspect
the number of lattice points from $\widehat L \cap P$ which differ
by translations by $e$, or by $f$. For an isolated bi-partition each
point outside of the parallelogram $\Sigma \subset P$
spanned by $e,f$ is translated  into $\Sigma$ by $e$ or by $f$, Fig. 5.
Second application of the same translation takes the point outside of
$P$. Hence the number of hinged pairs is equal to the number of points
outside of $\Sigma$, with the additional identification of centrally
symmetric points, which correspond to partitions equivalent by
a rotation. In this way we get the formula for the number of hinged pairs.
Let us note that special attention has to be paid to points in $P$
which have order two as elements of the group $\widehat L/L$. Such points
lying in the centers of the sides of the parallelogram $\Sigma$
do not give us a hinged pair because their translates produce a centrally
symmetric pair of points.

To obtain $s$ hinged Berg partitions we need to find a point in
$\widehat L \cap P$ such that $s-1$ of its translates by, say, $f$
are also in $P$. This implies
\[
(s-1)u < u+v, \ \ (s-1)q < p+q, \ \ s-2 < \min\{\frac{v}{u},\frac{p}{q}\},
\]
which gives us the required estimate.

In Case 1 the corner $-e'$ of the box $Q$ enters the box $P$ under
the translation by $f$.
To construct a matrix $C$ with many hinged partitions we
postulate that this corner
is translated in $s+1$ steps into the point $e+f+e'+f'$,
that is $e = sf -2e'-f'$.
(Note that
putting the endpoint of this string of points into the centrally symmetric
corner $e+f+e'$ would produce half as many hinged partitions, since
there would be equivalent pairs.)
Using the expression of $e,f$
in terms of $e',f'$ (elaborated in the proof of the previous theorem)
$e = (k-1)e' + l f', \ f = me'+ (n-1) f'$, we obtain
$k = sm-1, l = s(n-1)-1$. The requirement that $\det C = 1$ leaves
two free parameters $m(s+1)= n+1$ and
$C=
\left[ \begin{array}{rr} sm-1 & m
\\ s(s+1)m-2s-1 &  (s+1)m-1 \end{array} \right].
$
For such matrices we get
$s$ hinged Berg partitions. Putting $m=1$ we get the matrix above
for which the estimate is sharp. For $m \geq 2$ the estimate may be off by
$1$.
\end{Proof}$\square$

In particular we arrive at the conclusion that there is no universal
bound on the number of hinged Berg partitions.

Combining Proposition~\ref{wachlarz}, Theorem~\ref{repaut},
and the formula in Theorem~\ref{ilerozb} we can count the total number
of nonequivalent Berg partitions for a given automorphism, at least
for simple cutting sequences.

If the cutting sequence has the semiperiod $N$ with all the same elements,
say equal to $1$, and $\D$ generates its own centralizer
then the connectivity matrices are equal to
$C=
\left[ \begin{array}{rr} s & 1
\\ s(N-s)+1 &  N-s \end{array} \right],
$
for $s= 0,1,\dots,N-1$. We obtain
that  the number of all Berg partitions for $\D$ equals
$$\left\{\begin{array}{ll}
   \frac{1}{12}N(N^2+6N+5)\;:\;\;\;N\mbox{ is odd}\\&\\
   \frac{1}{12}N(N^2+6N+8)\;:\;\;\;N\mbox{ is even}.
   \end{array}\right.$$

Similarily if the cutting sequence has the period $N$ with $n_1$ ones
and $n_2$ zeroes, $n_1 \neq n_2$, then the total number of
nonequivalent Berg partitions is equal to
$$\left\{\begin{array}{ll}
   W(n_1,n_2)\;:\;\;\;n_1+n_2\mbox{ is even}\\&\\
   W(n_1,n_2)+\frac{1}{4}n_1\;:\;\;\;n_1+n_2\mbox{ is odd  and  $n_1$ is
even },
   \end{array}\right.$$
where $W(x,y)=\frac{1}{12}\left[(x^2+y^2)(xy+6)+6(x+y)(xy+1)+10xy\right].$

\section{Symmetries of bi-partitions }\label{symrozb}

Let us consider two irrational lines in the torus, as in Section~\ref{bipart}.
The construction of the bi-fan
$\F =\left\{ \{e_n,f_n\}|n \in \Z \right\} $,
and the respective sequence
of bi-partitions depends on the choice of the ``horizontal'' line.
When we exchange the roles of horizontal and vertical lines,
we obtain another bi-fan
$\widehat\F =\left\{ \{\hat e_n, \hat f_n\}|n \in \Z \right\} $,
where $\hat e_n = e_{-n}, \hat f_n = -f_{-n}$, $ n \in \Z$.
The bi-partitions associated with these bases are exactly the same but
their ordering is reversed. Let $\{s_n\}$ and $\{\hat s_n\}$ be the
respective cutting sequences. We have $\hat s_{n} = s_{-n-1}$.
It follows from the uniqueness of the cutting sequence that
\begin{proposition}\label{palindrom}
There is a toral automorphism (non-hyperbolic) which exchanges the
two irrational lines if and only if the sequences
$\{\hat s_n\}$ and $\{s_n\}$, or $\{1- s_n\}$, differ by a translation.
In the latter case the symmetry is of order 4, and in the former it
is an involution.
\end{proposition}
\begin{Proof}
By shifting the index appropriately we can assume that
one of the three possibilities takes place.
Firstly, $s_{n} = s_{-n-1} = \hat s_n, n\in \Z$, secondly
$s_{n} = s_{-n} = \hat s_{n-1}, n\in \Z$,
and thirdly $s_{n} = 1- s_{-n-1} = 1-\hat s_n, n\in \Z$.

In the first case the mapping which takes
$\{e_0,f_0\}$
into $\{\hat e_0, \hat f_0\} = \{e_0,-f_0\}$ exchanges the two lines,
and it has order 2.

In the second case
the mapping which takes
$\{e_0,f_0\}$
into $\{\hat e_{-1}, \hat f_{-1}\}$  is the desired
automorphism. We have either
$\{\hat e_{-1}, \hat f_{-1}\} = \{ e_{0}-f_{0},  -f_{0}\}$, if $s_{-1} =1$,
or $\{\hat e_{-1}, \hat f_{-1}\} = \{ e_{0},  e_{0}-f_{0}\}$, if $s_{-1}=0$.
In both cases we get an automorphism of order 2, nonconjugate to
the symmetry in the first case.

In the third case the mapping which takes
$\{e_0,f_0\}$
into $\{\hat f_0, \hat e_0\} = \{-f_0, e_0\}$ exchanges the two lines,
and it has order 4.
\end{Proof}$\square$

There are two distinct cases of the symmetry of order 2, just as there
are two nonconjugate elements of order 2 in $GL(2,\Z)$.
We will call the first type {\it simple order 2} symmetry,
and the second type {\it shift order 2 symmetry}. The geometric difference
is that a simple order 2 symmetry takes
the class of bi-partitions associated
with $\{e_0,f_0\}$ to itself, and pairs the other classes.
However no bi-partion in the invariant class is taken to itself.
With the shift order 2 symmetry no bi-partiton is taken to its translate,
all classes of bi-partitions are paired.
Note also that the  order of elements in a bi-partition is always reversed
by a symmetry of order 2, both simple and shift.

For the symmetry $\mathcal S$ of order 4 the class of bi-partitions
associated with $\{e_0,f_0\}$ is  mapped to itself,
and the other classes of bi-partitions are paired.
The invariant class contains exactly one
bi-partition which stays put under the symmetry, it happens when the
two fixed points of $\mathcal S$ coincide with the centers of the
rectangles.
For any bi-partition the order of elements is preserved under the symmetry.

In the case of symmetries of order 2, both simple and shift,
we can introduce coordinates so that the symmetry becomes the
Euclidean reflection in the diagonal.

In the case of a simple order 2 symmetry we get then that
$e_0 = (v,v) ,f_0 = (-u,u)$. The respective bi-partiton is taken to
itself with the exchange of the rectangles, wich are of sizes $u$  by $v$,
and $v$ by $u$. This bi-partition is automatically connected.
The vectors in the first quadrant
from the bi-fan $\mathcal F$ are symmetric with
respect to the diagonal, and the vectors in the second quadrant
are symmetric with respect to the anti-diagonal.

In the case of shift order 2 symmetry
we get either  $e_0 = (v,u+v), f_0 = (-u,u)$ or
$e_0 = (v,v), f_0 = (-u,u+v)$. Again the vectors in $\F$ are
symmetric with respect to the diagonal, and the anti-diagonal.

For the symmetry of order 4, we can choose coordinates so that
the symmetry becomes the rotation by
$\frac{\pi}{2}$. We get then $e_0 = (v,u), f_0 = (-u,v)$, and
the respective bi-partition contains two squares,
with sides $u$ and $v$ respectively. This bi-partition is
by necessity isolated. If the symmetry fixes
the center of one of the squares then it fixes the other center
as well, and both squares are mapped to themselves.

\section{Symmetries of toral automorphisms }

Let us now consider a hyperbolic toral automorphism $\D$ and
a symmetry $\mathcal S$ of the type discussed in Section~\ref{symrozb},
which exchanges the stable and unstable lines.
It is straightforward that $\mathcal S\circ \D \circ \mathcal S^{-1}
= \pm\D^{-1}$. Hence the automorphism is $\mathcal S$-reversible,
if $det ~\D  =1$. Moreover it follows that
$S$ takes Berg partitions of $\D$ into Berg partitons of $\pm\D^{-1}$.
Note that Berg partitions of $\pm \D$ and $\pm\D^{-1}$ have the same
shapes. However Berg partitions of $\D$ and $-\D$ are not the same.

Let the connectivity matrix for the automorphism $\D$ and
a Berg partition associated with
$\{e_{n},f_{n}\}$ be equal to $C_n$ for $n\in\Z$. Then
the connectivity matrix for $\D^{-1}$ and
the Berg partition associated with
$\{\hat e_{n}, \hat f_{n}\}$ is equal to $\sigma C_{-n}^T \sigma$.

For a $2\times 2$ matrix $A$ we define $A^J = \sigma A^T \sigma$,
i.e., $A^J$ is obtained from $A$ by the exchange of diagonal elements.
We say that $A$ is  {\it J-symmetric} if $A^J = A$, i.e., when $A$
has equal elements on the diagonal.

We get in the presence of a symmetry $S$ that the connectivity matrices
for two Berg partitions taken into one another by $S$ are
$C$ and $C^J$ in the case of symmetry of order 2, and $C$ and $C^T$ in the
case of symmetry of order 4.

For a simple order 2 symmetry $S$ there is class of Berg partitions
which are taken to their translates by $S$. For such a Berg partition
the connnectivity matrix is J-symmetric. For the symmetry of order 4
there is also a Berg partition which is taken to its translate, but
its connectivity matrix is symmetric.

Let us note that the presence of the simple order 2 symmetry
is characterized by the J-symmetry of a connectivity matrix,
and the presence of a symmetry of order 4 is characterized by
the symmetry of a connectivity matrix.

As shown in Section~\ref{symrozb} the special Berg partitions
which are taken into their translates by a symmetry have special shapes.
In the case of a simple order
2 symmetry by a choice of a flat metric in the torus they can be made into
two rectangles with sides $u$ by $v$ and $v$ by $u$.
In the case of a  symmetry of order 4 they can be made into
two squares.

There are six types of possible symmetries for hyperbolic automorphisms.

I. $\mathcal S$ is a simple order 2 symmetry,
the cutting sequence has the period $N=2k$.

In this case $\mathcal S$ takes $\{e_{n},f_{n}\}$
into $\{\hat e_{n},\hat f_{n}\}= \{e_{-n},-f_{-n}\}$ for every $n\in \Z$.
In particular it takes  $\{e_{k},f_{k}\}$
into $\{\hat e_{k},\hat f_{k}\}= \{e_{-k},-f_{-k}\}$.
The respective bi-partitions are equivalent
by the generator $\G$ of $\D$. It means that two equivalence classes
of Berg partitions,
those associated with $\{e_{0},f_{0}\}$ and $\{e_{k},f_{k}\}$
are fixed by $\mathcal S$ and the others are paired.

The symmetry $\mathcal S$ forces
$C_{n} = C_{-n}^J, n\in\Z$. In particular
the two matrices $C_0$ and $C_k$, and only them are J-symmetric.
Indeed, if $C_m$ is J-symmetric then $C_{-m} = C_m$, and that tell us
that $2m$ is a period.

The shortest period of the cutting sequence with this type of symmetry
has $N=6$, for instance $(110011)$.

Finally the sequence of $N=2k$ different connectivity matrices
for the different Berg partitions of $\D$ is
\[
C_0=C_0^J, C_1,\dots,C_{k-1},C_k = C_k^J, C_{k-1}^J,\dots,C_1^J.
\]

II. $\mathcal S$ is a simple order 2 symmetry,
the cutting sequence has the period $N=2k-1$.

Now there is also a shift order 2 symmetry equal to $\G\circ\mathcal S$.
More generally $\G^n\circ\mathcal S$ is a simple order 2 symmetry for
even $n$ and a shift order 2 symmetry for odd $n$.

The shortest period of the cutting sequence with this type of symmetry
has $N=5$, for instance $(11011)$. The sequence of $2k-1$ different
connectivity matrices is
\[
C_0=C_0^J, C_1,\dots,C_{k-1}, C_{k-1}^J,\dots,C_1^J.
\]

III. $\mathcal S$ is a shift order 2 symmetry, the cutting sequence
has the period $N=2k$

In this case $\G\circ\mathcal S$ is another
shift order 2 symmetry. There are no J-symmetric, or symmetric,
connectivity matrices.

The shortest period of the cutting sequence with this type of symmetry
is $N=4$, for instance $(1101)$. The sequence of $2k$ different
connectivity matrices is
\[
C_0, C_1,\dots,C_{k-1}, C_{k-1}^J,\dots,C_1^J,C_0^J.
\]

IV.  $\mathcal S$ is a symmetry of order 4, the cutting sequence
is periodic, and the period $N$ is by necessity even $N=2k$.

In this case $\G\circ\mathcal S$ is another
symmetry of order 4.
The shortest period of the cutting sequence with this type of symmetry
is $N=6$, for instance $(110100)$. The sequence of $2k$ different
connectivity matrices is
\[
C_0=C_0^T, C_1,\dots,C_{k-1}, C_k=C_k^T, C_{k-1}^T,\dots,C_1^T.
\]

V. $\mathcal S$ is a symmetry of order 4, the cutting sequence
has the semi-period $N=2k-1$.

In this case $\G\circ\mathcal S$ is a
shift order 2 symmetry.

The sequence of $N=2k-1$ different
connectivity matrices is
\[
C_0=C_0^T, C_1,\dots,C_{k-1}, C_{k-1}^J,\dots,C_1^J.
\]

If $N = 1$ and the semi-period is $(1)$
then in the basis $\{e_0,f_0\}$ we have $\G =
\left[ \begin{array}{rr} 0 & 1 \\ 1 &  1 \end{array} \right]$.
and $\mathcal S =
\left[ \begin{array}{rr} 0 & 1 \\ -1 &  0 \end{array} \right]$.

VI. $\mathcal S$ is a symmetry of order 4, the cutting sequence
has the semi-period $N=2k$.

In this case $\G\circ\mathcal S$ is a
simple order 2 symmetry.

The sequence of $N=2k$ different
connectivity matrices is
\[
C_0=C_0^T, C_1,\dots,C_{k-1},C_k=C_k^J, C_{k-1}^J,\dots,C_1^J.
\]

If $N =2$ and the semi-period is equal to $(11)$
then in the basis $\{e_0,f_0\}$ we have $\G =
\left[ \begin{array}{rr} 0 & 1 \\ 1 &  2 \end{array} \right]$.
and $\mathcal S =
\left[ \begin{array}{rr} 0 & 1 \\ -1 &  0 \end{array} \right]$.
If $\D =\G$ then the connectivity matrices are
\[
C_0 =\left[ \begin{array}{rr} 0 & 1 \\ 1 &  2 \end{array} \right],
C_1 =\left[ \begin{array}{rr} 1 & 1 \\ 2 &  1 \end{array} \right].
\]

Let us finally observe that the symmetries $\mathcal S$ and
$G\circ \mathcal S$ can be of one of three types. If we take into account
that they play the same role (i.e., their order is insignificant) then
we get six different ways of assigning the types, and all six are realized
as shown above.

At the same time, from the algebraic point of view, there are only
three different full symmetry groups for hyperbolic toral authomorphisms.
However, one of them is embedded in $Gl(2,\mathbb{Z})$ in three different
ways, another one in two ways, and the third one in one way only.
For details we refer to the paper of  \cite{B-R},
or the forthcoming  paper \cite{S-W}.


\begin{thebibliography}{M-R-R}



\bibitem[Ad]{Ad} R.L. Adler,
\emph{Symbolic Dynamics and Markov Partitions}
BAMS 35, (1998) 1-56


\bibitem[A-P]{A-P} R.L. Adler,  Palais,
\emph{Homeomorphic Conjugacy of Automorphisms on the Torus}
PAMS 16 (1965) 1222-1225



\bibitem[A-W]{A-W} R.L. Adler,B. Weiss,
\emph{Similarity of automorphisms of the torus}
Memoirs AMS 98 (1970)


\bibitem[A-K-K]{A-K-K} D. V. Anosov, A. V. Klimenko, G. Kolutsky, \emph{On the hyperbolic
    automorphisms of the 2-torus and their Markov partitions}, eprint arXiv:0810.5269

\bibitem[A-O]{A-O} H. Appelgate, H. Onishi,
\emph{Continued fractions and the conjugacy problem in ${\rm SL}_{2}(Z)$},
Comm. Algebra 9 (1981), no. 11, 1121--1130.


\bibitem[Ar]{Ar} D.Z. Arov,
\emph{Topological similitude of automorphisms of compact
commutative groups}
Uspehi Mat. Nauk 18 (1963)133-138. (Russian)


\bibitem[B-R]{B-R} M. Baake, A.G. Roberts, \emph{Reversing symmetry group of $Gl(2,\mathbb{Z})$ and $PGl(2,\mathbb{Z})$ matrices with connections to cat maps and trace maps}, J. Phys. A: Math. Gen. 30 (1997) 1549–1573.

\bibitem[Be]{Be} K. Berg,
\emph{On the conjugacy problem for K-systems}
Ph.D. Thesis, University of Minnesota, 1967.




\bibitem[Bo]{B} R. Bowen,
\emph{Markov partitions for Axiom A diffeomorphisms}
Amer. J. Math. 92 (1970), 725-747.

\bibitem[M]{M} A. Manning,
\emph{A Markov partition that reflects the geometry of a hyperbolic
toral automorphism}
TAMS 354 (2002), no. 7, 2849–2863 (electronic).





\bibitem[Se]{Se} C. Series,
\emph{The modular surface and continued fractions}
J. London Math. Soc. (2)  31  (1985),  no. 1, 69--80.

\bibitem[Si1]{Si1} Ya. G. Sinai,
\emph{Markov partitions and c-diffeomorphisms}
Funct. Anal. Applications 2 (1968), 61-82.

\bibitem[Si2]{Si2} Ya. G. Sinai,
\emph{Construction of Markov partitions}
Funct. Anal. Applications 2 (1968), 245-253. MR 40:3591


\bibitem[S-W]{S-W} A. Siemaszko, M.P. Wojtkowski,
\emph{The group centralizer and full symmetry group of toral endomorphisms}, in preparation.


\bibitem[Sn]{Sn} M. R. Snavely,
\emph{Markov partitions for the two-dimensional torus}
PAMS 113 (1991),  517--527.


\end{thebibliography}
\end{document}